\documentclass[12pt]{amsart}
\usepackage{fullpage,url}
\usepackage{amssymb}
\usepackage[all]{xy} 

\usepackage{longtable}

\DeclareFontEncoding{OT2}{}{} 
\newcommand{\textcyr}[1]{%
 {\fontencoding{OT2}\fontfamily{wncyr}\fontseries{m}\fontshape{n}\selectfont #1}}
\newcommand{\Sha}{{\mbox{\textcyr{Sh}}}}

\newcommand{\C}{{\mathbb C}}

\newcommand{\PP}{{\mathbb P}}
\newcommand{\Q}{{\mathbb Q}}

\newcommand{\Z}{{\mathbb Z}}

\newcommand{\Vbar}{{\overline{V}}}
\newcommand{\kbar}{{\overline{k}}}

\newcommand{\Pbar}{{\overline{P}}}

\newcommand{\calA}{{\mathcal A}}

\newcommand{\calL}{{\mathcal L}}

\newcommand{\calV}{{\mathcal V}}
\newcommand{\calW}{{\mathcal W}}

\DeclareMathOperator{\inv}{inv}

\DeclareMathOperator{\Aut}{Aut}
\DeclareMathOperator{\Gal}{Gal}

\DeclareMathOperator{\Br}{Br}

\DeclareMathOperator{\Sel}{Sel}

\DeclareMathOperator{\Div}{Div}
\DeclareMathOperator{\Pic}{Pic}

\DeclareMathOperator{\Spec}{Spec}

\DeclareMathOperator{\Jac}{Jac}

\newtheorem{theorem}{Theorem}[section]
\newtheorem{lemma}[theorem]{Lemma}

\newtheorem{proposition}[theorem]{Proposition}

\theoremstyle{definition}
\newtheorem{definition}[theorem]{Definition}

\newtheorem{remark}[theorem]{Remark}
\newtheorem{algorithm}[theorem]{Algorithm}

\begin{document}

\title[paper]{Tate-Shafarevich groups and K3 surfaces}
\author{Patrick Corn}
\address{Department of Mathematics, University of Georgia, 
	Athens, GA 30602-7403, USA}
\email{corn@math.uga.edu}
\urladdr{http://www.math.uga.edu/\~{}corn}
\date{\today}

\begin{abstract}

This paper explores a topic taken up recently by Logan and van Luijk in \cite{ronaldadam}--finding nontrivial $2$-torsion elements of the Tate-Shafarevich group of the Jacobian of a genus-$2$ curve by exhibiting Brauer-Manin obstructions to rational points on certain quotients of principal homogeneous spaces of the Jacobian, whose desingularizations are explicit K3 surfaces. The main difference between the methods used in this paper and those of Logan and van Luijk is that the obstructions are obtained here from explicitly constructed quaternion algebras, rather than elliptic fibrations. 
 
\end{abstract}

\maketitle

\section{Introduction}

Let $C$ be a curve of genus $2$ over a number field $k$, with Jacobian $J$. In an effort to describe the (finite) set $C(k)$, we are led to the study of $J(k)$. To determine its rank, we refer to a well-known exact sequence 
\begin{equation}\label{exact}
0 \to J(k)/2J(k) \to \Sel^{(2)}(k,J) \to \Sha(k,J)[2] \to 0
\end{equation}
where the middle group is effectively computable. See \cite{stoll} for a comprehensive description of the computation, which has been implemented in the computer algebra system MAGMA (\cite{magma}).

So computing the group on the left is more or less equivalent to computing the rather mysterious group on the right. In this paper, we find examples of curves $C$ over $\Q$ such that $\Sha(\Q,J)[2]$ is nonzero, by finding explicit elements of this group. Such elements can be represented by $2$-coverings $X$ of $J$ which have points everywhere locally but no $k$-points. The strategy, as in \cite{ronaldadam}, is to prove that the Hasse principle fails for $X$ by exhibiting a Brauer-Manin obstruction to rational points on the desingularization of the quotient $X/\iota$, where $\iota$ is the involution corresponding to multiplication by $-1$. This desingularization is a K3 surface which can be given explicitly as the smooth complete intersection of $3$ quadrics in $\PP^5$.

The result is the following theorem.

\begin{theorem} \label{main} Let $S$ be the set of primes splitting completely in a certain finite extension $K/\Q$ (this extension is given explicitly in the statement of Proposition \ref{mainthm}). For all $n$ equal to the product of primes in $S$, the genus-$2$ curve 
\[
C_n \colon y^2 = n(x^2-5x+1)(x^3-7x+10)(x+1)
\]
satisfies $\Sha(\Q,\Jac(C))[2] \ne 0$. \end{theorem}

The curve with $n=1$ was originally obtained by a computer search over products of polynomials of degree $2,3,1$ with small coefficients.

I thank Ronald van Luijk for introducing me to this topic and for many helpful conversations. I would also like to thank Bjorn Poonen for several enlightening comments, and in particular for the method of computing problematic smooth places outlined in the proof of the main theorem. 

\section{The Brauer-Manin obstruction}

\subsection{Generalities} First we briefly review the Brauer-Manin obstruction to the Hasse principle. Let $V$ be a smooth proper $k$-variety, $k$ a number field. Then the map $V({\mathbb A}_k) \to \prod_v V(k_v)$ is a bijection (\cite{skorobogatov}, pp. 98-99). We will suppose that this set $V({\mathbb A}_k)$ is nonempty. 

For any scheme $V$ we can define the {\em Brauer group} $\Br V = H^2(V_{\rm et}, {\mathbb G}_m)$. For an element $\calA \in \Br V$, define the set
\[
V({\mathbb A}_k)^\calA = \{ (P_v) \in V({\mathbb A}_k) : \sum_v \inv_v \, \calA(P_v) = 0 \},
\]
where the sum is over all places $v$ of $k$, and define 
\begin{equation}\label{brint}
V({\mathbb A}_k)^{\rm Br} = \bigcap_{\calA \in \Br V} V({\mathbb A}_k)^\calA.
\end{equation}

Here $\inv_v \colon \Br k_v \to \Q/\Z$ is an isomorphism if $v$ is nonarchimedian, or the injection $\frac12\Z/Z \to \Q/\Z$ if $k_v = \mathbb R$, or the zero map if $k_v = \mathbb C$. Class field theory shows that 
\[
V(k) \subseteq V({\mathbb A}_k)^{\rm Br} \subseteq V({\mathbb A}_k),
\]
so we say that $V$ has a Brauer-Manin obstruction to the Hasse principle if $V({\mathbb A}_k)^{\rm Br}$ is empty, so that $V(k)$ is as well.

For many classes of varieties, it is believed that the Brauer-Manin obstruction to the Hasse principle is ``the only one"; that is, if $V({\mathbb A}_k)^{\rm Br}$ is nonempty, then so is $V(k)$. It is not known (even conjecturally) whether or not the Brauer-Manin obstruction to the Hasse principle is the only one for K3 surfaces.

\subsection{The key isomorphism}

The map $V \to \Spec k$ induces a map $\Br k \to \Br V$, which is injective if $V({\mathbb A}_k)$ is nonempty. Elements in the image of this map are called {\em constant algebras}. Two elements of $\Br V$ which differ by a constant algebra cut out the same subset of $V({\mathbb A}_k)$, so the intersection (\ref{brint}) defining $V({\mathbb A}_k)^{\rm Br}$ need only be taken over a set of representatives of $\dfrac{\Br V}{\Br k}$.

\begin{proposition} For $V$ a smooth projective geometrically integral $k$-variety, $k$ a number field, with $V({\mathbb A}_k) \ne \emptyset$, there is an isomorphism
\begin{equation}\label{briso}
\frac{\Br_1 V}{\Br k} \to H^1(k,\Pic \Vbar)
\end{equation}
where $\Vbar = V \times_k \kbar$ and $\Br_1 V = \ker(\Br V \to \Br \Vbar)$ is the ``algebraic part" of the Brauer group. \end{proposition}

{\em Proof:} This is a standard consequence of the Hochschild-Serre spectral sequence; see for instance \cite{thesis}, Proposition 1.3.7. $\Box$

\smallskip

It is, unfortunately, very difficult in general to compute the inverse of the isomorphism (\ref{briso}) explicitly. The crux of the computation is an explicit use of the fact (due to Tate) that $H^3(k,\kbar^*) = 0$; that is, one must express an arbitrary $3$-cocycle with values in $\kbar^*$ as a coboundary. In the next section, we discuss one way around this problem.

\subsection{Quaternion algebras in $\Br_1(V)$} One common way of constructing explicit elements of $\Br_1(V)$ is as follows: 

\begin{definition} For $c \in k^*$ and $g \in k(V)^*$, the quaternion algebra $(c,g)$ is a four-dimensional central simple $k(V)$-algebra with $k(V)$-basis $1,i,j,ij$, satisfying $i^2 = c$, $j^2 = g$, and $ij = -ji$. \end{definition}

Of course, quaternion algebras are quite general objects; the reason we study the special quaternion algebras defined above is the following standard lemma:

\begin{lemma} For $c \in k^*$ and $g \in k(V)^*$, a quaternion algebra $(c,g) \in \Br k(V)$ is in the image of the map $\Br V \to \Br k(V)$ if and only if $\text{div}(g) = D + \sigma D$, where $D$ is a divisor defined over $k(\sqrt{c})$ and $\sigma$ is the nontrivial element of $\Gal(k(\sqrt{c})/k)$. It is a constant algebra if and only if $\text{div}(g) = D' + \sigma D'$, where $D'$ is a {\em principal} divisor defined over $k(\sqrt{c})$. \end{lemma}

{\em Proof:} See \cite{thesis}, Proposition 2.2.3 or \cite{bright}, Proposition 4.17. $\Box$

\smallskip

Note that the quaternion algebra $(c,g)$ will always split over the quadratic extension $k(\sqrt{c})$, so its image in $H^1(k,\Pic \Vbar)$ will restrict to $0$ in $H^1(k(\sqrt{c}),\Pic \Vbar)$. There is also a well-known formula for the local invariant of such a quaternion algebra in $\Br V$: for any point $P_v \in V(k_v)$, we have that 
\[
\inv_v (c,g)(P_v) = [c,g(P_v)]_v,
\]
where $[a,b]_v$ is the Hilbert symbol of $a,b \in k_v^*$ expressed as an element of $\dfrac12\Z/\Z$. So $[a,b]_v = 0$ if and only if $x^2-ay^2-bz^2$ represents $0$ in $k_v$; otherwise it equals $1/2$. 

\begin{algorithm}\label{alg} Given a class of varieties over $k$, here is how we look for varieties $V$ in that class with quaternion algebras of the above type generating nonconstant elements in $\Br_1(V)$:

\begin{enumerate}
\item Find a $G_k$-invariant generating set $\Gamma$ for $\Pic \Vbar$ (possibly a subgroup of $\Pic \Vbar$ will work as well; see the comments). 
\item The action of the Galois group on $\Gamma$ induces a map $G_k \to \Aut(\Gamma)$. By inflation-restriction, $H^1(k,\Pic \Vbar)$ will be isomorphic to $H^1(H,\Pic \Vbar)$, where $H$ is the image of $G_k$ inside $\Aut(\Gamma)$. List the cohomology groups $H^1(H,\Pic \Vbar)$ for every subgroup $H$ of $\Aut(\Gamma)$; these are the possibilities for $H^1(k,\Pic \Vbar)$.
\item Search the subgroups of $\Aut(\Gamma)$ to find $H$ such that $H^1(H,\Pic \Vbar) = \Z/2$ but $H^1(H',\Pic \Vbar) = 0$ for some subgroup $H' \subset H$ of index $2$. 
\item Give conditions on $V$ such that the image of $G_k$ in $\Aut(\Gamma)$ equals $H$.
\item Let $L$ be the fixed field of $H'$. Use the equality $H^1(H,\Pic \Vbar) = \Z/2$ to find a nonzero divisor class $d \in \Pic V_L$ such that $d + \sigma d = 0$ in $\Pic V_L$, $\sigma$ the nontrivial element of $\Gal(L/k)$.
\item Given $d$ and $L$, find a divisor $D$ defined over $L$ whose class is $d$. (Such a divisor is guaranteed to exist whenever $V$ has points everywhere locally: see \cite{thesis}, Proposition 1.3.4.) Then our quaternion algebra is $(c,g)$, where $L = k(\sqrt{c})$ and $g$ is a function whose divisor is $D + \sigma D$.
\end{enumerate}

\end{algorithm}

{\em Comments on the algorithm:} Carrying out step (1) requires that we know quite a lot about the geometry of $\Vbar$. Even for general K3 surfaces, this is too hard. Some previous attempts to carry out this algorithm explicitly have restricted themselves to Del Pezzo surfaces (e.g. \cite{thesis}, \cite{kreschtschinkel2004}, \cite{cornlms}), or special K3 surfaces such as diagonal quartic hypersurfaces (\cite{bright}) or Kummer surfaces (\cite{argentinthesis}). As we will see, the K3 surfaces we examine in this paper come equipped with a $G_k$-invariant set $\Gamma$ of $32$ lines generating a subgroup of $\Pic \Vbar$ which is free of rank $17$.

In Step (3), we usually search for a subgroup $H$ maximal with respect to the given properties, so that we can study the broadest possible class of varieties $V$ admitting a nonconstant quaternion algebra of the above type in $\Br_1 V$.

We carry out Step (5) using the following procedure: let $\Pbar = \Pic \Vbar$. Then there is an isomorphism 
\[
\frac{(\Pbar/2\Pbar)^{G_k}}{\Pbar^{G_k}/2\Pbar^{G_k}} \to H^1(k,\Pbar)[2]
\]
sending (the class of) a divisor class $e$ on the left to the cocycle $c_{\tau} = \frac12(\tau e - e)$ (cf. \cite{cornlms}, Lemma 3.1). Since we have computed $H^1(k,\Pbar)[2]$, we can find a nontrivial $c_{\tau}$ and hence a nontrivial $e$. If this cocycle trivializes upon restriction to $G_L$, then we can find a $e$ on the left defined over $L$, and then the divisor class $d = \frac12(\sigma e - e)$ satisfies $d + \sigma d = 0$ in $\Pic V_L$.

Moreover, in this case the class in $\Br_1(V)/\Br k$ of the quaternion algebra $(c,g)$ we obtain at the end of the algorithm corresponds via the isomorphism \ref{briso} to the class of the cocycle $c_{\tau}$, so if we know that this cocycle is not a coboundary, we are guaranteed that the quaternion algebra we have found is nonconstant.

In practice, we usually write down conditions on $V$ in Step (4) in such a way that {\em generically} $H^1(k,\Pbar)[2]$ is nontrivial, because the image of $G_k \to \Aut(\Gamma)$ is the subgroup $H$ we found in Step (3); but if the image of $G_k \to \Aut(\Gamma)$ is strictly contained in $H$, it might happen that the quaternion algebra $(c,g)$ we construct is constant.

Similarly, if $\Gamma$ generates a subgroup $Q$ of $\Pbar$, this algorithm constructs a quaternion algebra whose corresponding cocycle in $H^1(k,Q)$ is not a coboundary. If $Q$ is a proper subgroup of $\Pbar$, the restriction of this quaternion algebra to $H^1(k,\Pbar)$ might be zero. So in this case (and the problem case in the above paragraph), we will get a bona fide element of $\Br V$ from the above algorithm, but it will not furnish a Brauer-Manin obstruction to rational points on $V$.

\section{Genus-$2$ curves and K3 surfaces: a review of the theory}

The material in this section is taken largely from \cite{ronaldadam}. 

\subsection{Construction of the K3 surface} Let $C$ be a genus-$2$ curve over a number field $k$ and $J$ its Jacobian. We wish to associate, to an element $\delta \in \Sel^{(2)}(k,J)$, a K3 surface $V_{\delta}$ with points everywhere locally, such that if $\delta$ is in the image of $J(k)/2J(k)$, then $V_{\delta}$ has a rational point.

\begin{definition} For $f(x) \in k[x]$ of degree 6, set $A_f = k[x]/(f(x))$. \end{definition}

\begin{definition} For any $\delta \in A_f^*$, define
\[
\calV_{f,\delta} = \{ q \in A_f \colon \delta q^2 \equiv \text{quadratic mod $f$} \}.
\]
\end{definition}

If we write $q(x) = \sum_{i=0}^5 a_i x^i$, then $q$ lies in $\calV_{f,\delta}$ if and only if the coefficients of $x^3$, $x^4$, and $x^5$ vanish. Considered as polynomials in the $a_i$, these coefficients $C_3, C_4, C_5$ are homogeneous of degree $2$. 

\begin{definition} For $\delta \in A_f^*$, define $V_{\delta}$ to be the variety in $\PP^5$ consisting of points $(a_0 \colon \cdots \colon a_5)$ such that $\sum_{i=0}^5 a_i x^i \in \calV_{f,\delta}$. One shows (\cite{ronaldadam}, Proposition 2.1.2) that $V_{\delta}$ is a smooth complete intersection of the three quadrics $C_3, C_4, C_5$, which makes it a K3 surface of degree $8$. \end{definition}

Let $H_f$ be the kernel of the norm map $A_f^*/(A_f^*)^2 k^* \to k^*/(k^*)^2$ (which is well-defined because the degree of $f$ is even). Then lemma 5.1 of \cite{stoll} describes a map $\Delta_k \colon J(k)/2J(k) \to H_f$ whose kernel has order $1$ or $2$; as usual, it is induced from the homomorphism $\Div_{\perp}^0(C)(k) \to A_f^*$ defined by $P \mapsto x(P)-x$, where $\Div_{\perp}^0(C)$ consists of divisors whose support is disjoint from that of $\text{div}(y)$. In addition, the lemma gives necessary and sufficient conditions under which the kernel will have order $1$. It suffices, for instance, for $f$ to have a factor of odd degree. Cf. Theorem 13.2 of \cite{poonenschaefer}; following section 5 of \cite{stoll}, we say that ``$k$ satisfies condition $(\ddagger)$" (thinking of $f$ as fixed) if the kernel has order $1$. 

If $k$ satisfies condition $(\ddagger)$, we can identify $\Sel^{(2)}(k,J)$ with the set of elements of $H_f$ whose images in $H_f \otimes k_v$ lie in the image of $\Delta_{k_v}$ for all places $v$ of $k$.\footnote{In general, this set is known as the ``fake $2$-Selmer group," but it is isomorphic with the $2$-Selmer group if $k$ satisfies condition $(\ddagger)$--see \cite{stoll}, section 5.}  

Henceforth, we will assume that $k$ satisfies $(\ddagger)$, and we will think of $\Sel^{(2)}(k,J)$ as the subset of $H_f$ described in the above paragraph. The goal is then to compute that set and to look for elements in it which are not in the image of $\Delta_k$.

This is where the set $\calV_{f,\delta}$ comes in: Proposition 3.2.6 of \cite{ronaldadam} shows that if $\delta \in \Sel^{(2)}(k,J)$ is in the image of $\Delta_k$, then there is a polynomial $q \in \calV_{f,\delta}$. For the sake of concreteness, we sketch the proof: $\delta$ will be congruent mod $(A_f^*)^2 k^*$ to $(x(P_1)-x)(x(P_2)-x)$ for some pair of points $P_1,P_2 \in C(\kbar)$ which are defined and conjugate over a quadratic extension of $k$. We restate this fact:

\begin{proposition}\label{points} Let $C$ be a genus-$2$ curve given by $y^2 = f(x)$, $f(x) \in k[x]$ of degree $6$, and let $J$ be its Jacobian. Suppose that $k$ satisfies $(\ddagger)$. If $\delta \in \Sel^{(2)}(k,J)$ is in the image of the map $J(k)/2J(k) \to \Sel^{(2)}(k,J)$ from (\ref{exact}), the K3 surface $V_{\delta}$ has a rational point. \end{proposition} 

We can view this result from another perspective as well: elements $\delta \in \Sel^{(2)}(k,J)$ correspond to $2$-coverings $X_{\delta}$ of $J$, which are twists of $J$ equipped with a map $\pi \colon X \to J$ defined over $k$ which is a twist of multiplication by $2$ on $J_{\kbar}$. Such $2$-coverings inherit an involution defined over $k$ descending from multiplication by $-1$ on $J$. Then the surface $V_{\delta}$ is the minimal nonsingular model of the quotient of $X_{\delta}$ by this involution. (It is a twist of the desingularized Kummer surface of $J$; cf. Chapter 16 of \cite{prolegomena}.) 

The element $\delta$ comes from $J(k)/2J(k)$ if and only if $X_{\delta}$ is the trivial twist of $J$, which happens if and only if $X_{\delta}(k)$ is nonempty. But since $\delta \in \Sel^{(2)}(k,J)$, we must have that $X_{\delta}(k_v)$ is nonempty for all places $v$ of $k$. Thus $X_{\delta}$ is a counterexample to the Hasse principle. Certainly if $X_{\delta}$ has points in some field, then so does $V_{\delta}$; so $V_{\delta}$ has points everywhere locally, and $V_{\delta}(k) = \emptyset$ will imply that $X_{\delta}(k) = \emptyset$ and hence that $\delta$ maps to a nontrivial element of $\Sha(k,J)[2]$.

\begin{remark} It may be true that $X_{\delta}(k)$ is empty while $V_{\delta}(k)$ is nonempty; since $X_{\delta}$ is really the surface in whose rational points we are interested, one might hope to work directly with it instead of $V_{\delta}$. Of course, the reason we do not do this is that the usual explicit description of $X_{\delta}$ is as an intersection of $72$ quadrics in $\PP^{15}$ (just as it is for the Jacobian itself--see \cite{prolegomena}, Chapter 2 for the construction). \end{remark}

\begin{remark}\label{dp4} If $f(x)$ has a $k$-rational root, we can find a model $C'$ for $C$ of the form $w^2 = g(u)$, deg $g = 5$. If we carry through the above constructions in $A_g$ instead of $A_f$, we get a smooth complete intersection of two quadrics in $\PP^4$, which is a Del Pezzo surface $W_{\beta}$ of degree $4$, where $\beta$ is the element of $\Sel^{(2)}(k,\Jac C')$ corresponding to $\delta$ (see \cite{bruinflynn} and \cite{logandp4}). We will see later that $V_{\delta}$ is a double cover of $W_{\beta}$. Then, just as in the previous remark, it may be true that $W_{\beta}(k)$ is nonempty while $V_{\delta}(k)$ is empty; in fact, this happens for the curve given in Theorem \ref{main} (assuming that the Brauer-Manin obstruction to the Hasse principle is the only one for Del Pezzo surfaces).
\end{remark}

\subsection{The 32 lines and $\Pic {\overline {V_{\delta}}}$} Here we review the construction of the $32$ lines on $\Pic {\overline {V_{\delta}}}$ and analyze the structure of $\Pic {\overline {V_{\delta}}}$ as a $G_k$-module. The ideas are taken from \cite{ronaldadam}, but the notation will be different.

For $\delta \in A_f^*$, let $r_1, \ldots, r_6$ be the roots of $f$ in $\kbar$. Fix a choice of square roots $z_i$ of $\delta(r_i)$ in $\kbar$.

For an element $s = (s_1, \ldots, s_6) \in (\Z/2)^6$, define $\gamma_s$ to be the unique degree-$5$ polynomial satisfying $\gamma_s(r_i) = \dfrac{(-1)^{s_i}}{z_i}$ for $1 \le i \le 6$. Now define
\[
\calL_s = \{ \gamma_s(x)(tx+u) \colon t,u \in \kbar \}
\]
and notice that $\delta(x) q(x)^2 \equiv (tx+u)^2$ mod $f$, for any $q(x) \in \calL_s$. The projectivization of $\calL_s$ is a line $L_s$ on $\overline {V_{\delta}}$. Notice that if $s + s' = (1,1,1,1,1,1)$, we have that $\gamma_s = -\gamma_{s'}$, so $L_s = L_{s'}$. Thus we have $32$ lines $L_s$ indexed by elements of $(\Z/2)^6/\langle (1,1,1,1,1,1) \rangle$.

It is not hard to determine the intersection pairing as applied to the subgroup generated by these lines in $\Pic \overline {V_{\delta}}$; one obtains 
\[
L_{s_1} \cdot L_{s_2} = \begin{cases} -2 & \text{if $s_1$ and $s_2$ differ in $0$ or $6$ places} \\
1 & \text{if $s_1$ and $s_2$ differ in $1$ or $5$ places} \\
0 & \text{otherwise} \end{cases}
\]

The first line follows by the adjunction formula for K3 surfaces, and the second and third lines by direct calculations. From these intersection numbers we obtain

\begin{proposition} (\cite{ronaldadam}, Proposition 2.1.21) The classes of the $32$ lines on $\overline{V_{\delta}}$ generate a subgroup of $\Pic \overline{V_{\delta}}$ isomorphic to $\Z^{17}$. \end{proposition}

\begin{remark} In the generic case this subgroup equals all of $\Pic \overline{V_{\delta}}$ (Proposition 2.1.30 of \cite{ronaldadam}). As noted above, even if $V_{\delta}$ is not generic, we can still use the subgroup generated by the classes of the lines to give us an element of $H^1(k,\Pic \overline{V_{\delta}})$ via restriction. 
\end{remark}

\section{The algorithm for $V_{\delta}$}

Viewed in terms of Algorithm \ref{alg}, the end of the previous section carries out Step (1) for the class of varieties of the form $V_{\delta}$; then $\Gamma$ is in our case the set of $32$ lines on $\overline{V_{\delta}}$ defined earlier. We proceed to carry out the remaining steps of the algorithm.

Let $G_{\Gamma}$ be the group $(\Z/2)^6/\langle (1,1,1,1,1,1) \rangle$ indexing the $32$ lines. Then, by a straightforward analysis of the intersection pairing on $\Gamma$ (which is Proposition 2.2.11 of \cite{ronaldadam}), $\Aut(\Gamma)$ is isomorphic to $G_{\Gamma} \rtimes S_6$, where the symmetric group acts on $G_{\Gamma}$ by permuting indices (which corresponds to permuting the square roots $z_i$ of $\delta(r_i)$), and $G_{\Gamma}$ acts on itself by addition.

For $\delta \in \Sel^{(2)}(k,J)$, the image of $G_k \to \Aut(\Gamma)$ must lie in a subgroup of index $2$ inside $\Aut(\Gamma)$, because the norm of $\delta$ is required to be a square in $k$. Indeed, this index-$2$ subgroup is the semi-direct product of $S_6$ with the index-$2$ subgroup of $G_{\Gamma}$ of elements whose sum is $0$.

This is a group $G$ of order $11520$ which acts on $\Z^{17}$ in a prescribed way. MAGMA can enumerate its subgroups $H$ and compute $H^1(H,\Z^{17})$ for each one. This is Step (2).

We look in Step (3) for maximal subgroups $H$ such that $H^1(H,\Z^{17}) = \Z/2$ and $H^1(H',\Z^{17}) = 0$ for some index-$2$ subgroup $H' \subset H$. MAGMA finds two conjugacy classes of such subgroups. One class consists of subgroups of order $96$; the other consists of subgroups of order $128$. Here we exhibit a Brauer-Manin obstruction coming from the first conjugacy class; presumably we could use the same techniques to try to find one coming from the second class. 

\subsection{The subgroup $H_{96}$ and the shape of $f(x)$} Consider polynomials $f(x)$ of the shape
\[
f(x) = f_2(x)f_3(x)(x-r_6)
\]
where $f_i$ is an irreducible polynomial of degree $i$. Let $C$ be the genus-$2$ curve $y^2 = f(x)$. Pick the ordering of the roots of $f$ that lists the two roots of $f_2$, then the three roots of $f_3$, then $r_6$. Suppose also that $\delta \in \Sel^{(2)}(k,J)$ satisfies the condition that $\delta(r_1)$ is a square in $\Q(r_1)$, where $r_1$ is a root of $f_2(x)$. Then the image of $G_k \to \Aut(\Gamma)$ will sit inside the subgroup of $G_{\Gamma} \rtimes S_6$ generated by the elements
\[
(0,0,1,0,0,1), (0,0,1,1,1,1), (12), (345), (34)
\]
where the first two elements are in $G_{\Gamma}$ and the last three elements are permutations in $S_6$.

It is not hard to check that this is a subgroup $H_{96}$ of order $96$. If we define $H_{48}$ to be the index-$2$ subgroup consisting of elements which leave $z_6 = \sqrt{\delta(r_6)}$ unchanged, then MAGMA computes that $H^1(H_{96},\Z^{17}) = \Z/2$ and $H^1(H_{48},\Z^{17}) = 0$.

Now let us see where the nontrivial element of $H^1(H_{96},\Z^{17})$ comes from: set 
\[
d = \left( \ell_{(0,0,0,1,1,0)} + \ell_{(0,0,0,1,1,1)} \right) - \left( \ell_{(0,0,1,0,0,0)} + \ell_{(0,0,1,0,0,1)} \right),
\]
where $\ell_s$ is the divisor class of $L_s$. One computes that $\sigma d = d$ for all $d \in H_{48}$, and $\sigma d = -d$ for all $\sigma \in H_{96} \setminus H_{48}$. So $d$ is the divisor class we referred to in Step (5) of Algorithm \ref{alg}.

Now the fixed field of the intersection of the image with $H_{48}$ will be $L = k(\sqrt{\delta(r_6)})$. The remaining step in the algorithm is to find a divisor $D \in \Div (V_{\delta})_L$ whose divisor class is $d$. Of course we will have to use the fact that $V_{\delta}$ has points everywhere locally. The easiest way to solve this computational problem, generally speaking, is to reduce it to solving a certain norm equation which we are guaranteed has a solution by the Hasse Norm Theorem. We now show how this can be accomplished for $V_{\delta}$.

\subsection{The divisor $D$ and the quaternion algebra} The strategy is to take a divisor $E$ in the class of $d$ defined over a higher-degree extension of $k$, and then to subtract the divisor of a judiciously chosen rational function to $E$ in order to obtain a divisor defined over $L$. That is, we want
\[
\tau(E-(h)) = E-(h) \ \text{for all $\tau \in G_L$.}
\]
Here we begin with 
\[
E = \left( L_{(0,0,0,1,1,0)} + L_{(0,0,0,1,1,1)} \right) - \left( L_{(0,0,1,0,0,0)} + L_{(0,0,1,0,0,1)} \right),
\]

The stabilizer $H_{12}$ of the divisor $E$ (not its class) in $H_{96}$ is generated by $(0,0,1,1,0,0) \cdot (34), (12), (45)$. It has order $12$, and index $4$ in $H_{48}$. We wish to find a function $(h)$ such that $E-(h)$ is fixed by $H_{48}$. That is, for all $\tau \in H_{48}$, we want $E-\tau E = \text{div}(h/\tau h)$. 

There are four left cosets of $H_{12}$ in $H_{48}$, and we first identify what the divisor of $h/\tau h$ must be for $\tau$ in each coset. For ease of notation, we identify an element $s \in G_{\Gamma}$ with the binary number $\sum_{i=1}^6 2^{6-i} s_i$. So in this notation
\[
E = L_6 + L_7 - (L_8 + L_9).
\]

So we get
\begin{align*}
\tau \in H_{12} &\Rightarrow \ \text{div}(h/\tau h) = 0 \\
\tau \in (0,0,1,1,0,0) H_{12} &\Rightarrow \ \text{div}(h/\tau h) = L_4 + L_5 + L_6 + L_7 - (L_8 + L_9 + L_{10} + L_{11}) \\
\tau \in (0,0,0,1,1,0) H_{12} &\Rightarrow \ \text{div}(h/\tau h) = L_6 + L_7 + L_{14} + L_{15} - (L_0 + L_1 + L_8 + L_9) \\ 
\tau \in (0,0,1,0,1,0) H_{12} &\Rightarrow \ \text{div}(h/\tau h) = L_2 + L_3 + L_6 + L_7 - (L_8 + L_9 + L_{12} + L_{13}) 
\end{align*}

\begin{proposition} Fix an index set $I = i_1, i_2, i_3 \in \{ 1,2,3,4,5,6 \}$. Fix a sequence $B = (b_1,b_2,b_3) \in \Z/2$. Let $E_{I,B}$ be the sum of the eight lines $L_s$, $s \in G_{\Gamma}$, where $s_{i_j} = b_j$. Then $E_{I,B}$ is a hyperplane section. \end{proposition}

{\em Proof:} This is a restatement of Lemma 2.1.24 in \cite{ronaldadam}. $\Box$ 

\smallskip

The point of this proposition is that we can write the above divisors as differences of hyperplane sections. For $3 \le i \le 6$, define $p_i$ to be a linear polynomial cutting out $E_{\{1,2,i\},\{0,0,0\}}$; and define $q_i$ to be a linear polynomial cutting out $E_{\{1,2,i\},\{0,0,1\}}$. Then we see that
\begin{align*}
L_4 + L_5 + L_6 + L_7 - (L_8 + L_9 + L_{10} + L_{11}) &= \text{div} \left( \frac{p_3}{p_4} \right) \\
L_6 + L_7 + L_{14} + L_{15} - (L_0 + L_1 + L_8 + L_9) &= \text{div} \left( \frac{q_5}{p_4} \right)
\\
L_2 + L_3 + L_6 + L_7 - (L_8 + L_9 + L_{12} + L_{13}) &= \text{div} \left( \frac{p_3}{p_5} \right)
\end{align*}

The proof of Hilbert Theorem 90 suggests the following choice of $h$:
\[
h = 1 + \frac{p_3}{p_4} + \frac{q_5}{p_4} + \frac{p_3}{p_5}.
\]

Subject to some conditions on the $p_i$ and $q_i$, which have so far only been defined up to scalar multiples, we will show that this $h$ gives us the rational function we want. Here we outline the conditions we will need.

\smallskip

{\bf Condition 1:} We require $p_i$ and $q_i$ to have coefficients in $k(r_1,z_i,z_6)$, and to be defined so that $(0,0,1,1,1,1)p_i = q_i$. This condition will be immediate from the construction of the $p_i$ and $q_i$ we will outline below.

\smallskip

{\bf Condition 2:} We need that $\sigma p_i = p_j$ and $\sigma q_i = q_j$ if $\sigma \in H_{96} \cap S_6$ sends $i$ to $j$, and so that 
\[
\tau p_i = \begin{cases} p_i & \text{if $\tau \in H_{96} \cap G_{\Gamma}$ has a $0$ in the $i$th spot} \\
q_i & \text{if $\tau \in H_{96} \cap G_{\Gamma}$ has a $1$ in the $i$th spot} \end{cases}
\]
This condition is easy to satisfy, as we will shortly see.

\smallskip

{\bf Condition 3:} We also want that $p_i q_i = p_j q_j$ for $3 \le i,j \le 5$. This condition is more difficult to satisfy; we will need to use the fact that $V_{\delta}$ has points everywhere locally.

\smallskip

\begin{proposition} If $p_i$ and $q_i$ satisfy Conditions 1,2,3, then the divisor $E - (h)$ is defined over the field $L = k(\sqrt{\delta(r_6)})$. \end{proposition}

{\em Proof:} By condition 1, the rational function $h$ is defined over the field $k(r_1,z_3,z_4,z_5)$. Note that $z_6$ is automatically contained in this field extension because of the requirement that the product of the $z_i$ is in $k$ (and the requirement that $z_1 \in k(r_1)$). The Galois group of this field extension is a subgroup of $H_{96}$. So we must only show that $E-\tau E = {\rm div}(h/\tau h)$ for every $\tau \in H_{48}$.

First we show that $\tau h = h$ for $\tau \in H_{12}$. Clearly this is the case for $\tau = (12)$. Now 
\[
(45) h = 1 + \frac{p_3}{p_5} + \frac{q_4}{p_5} + \frac{p_3}{p_4} = 1 + \frac{p_3}{p_5} + \frac{p_4}{q_5} + \frac{p_3}{p_4} = h
\]
and
\begin{align*}
((0,0,1,1,0,0) \cdot (34)) h &= (0,0,1,1,0,0) \left( 1 + \frac{p_4}{p_3} + \frac{q_5}{p_3} + \frac{p_4}{p_5} \right) \\
&= 1 + \frac{q_4}{q_3} + \frac{q_5}{q_3} + \frac{q_4}{p_5} \\
&= 1 + \frac{p_3}{p_4} + \frac{p_3}{p_5} + \frac{p_4}{q_5} = h.
\end{align*}

So $\text{div}(h/\tau h)$ depends only on the left coset of $H_{12}$ to which $\tau$ belongs. In the computations that follow, it will help to note that
\[
h = \frac{p_4p_5 + p_3p_5 + p_3p_4 + p_5q_5}{p_4p_5}
\]
and to notice that the numerator is invariant under permutations of the coordinates $3,4,5$, by Condition 3. So we are reduced to three computations:
\begin{align*}
\frac{h}{(34) h} &= \frac{p_3p_5}{p_4p_5} = \frac{p_3}{p_4} \\
\frac{h}{(0,0,0,1,1,0) h} &= \frac{(p_4p_5 + p_3p_5 + p_3p_4 + p_5q_5)q_4q_5}{(q_4q_5 + p_3q_5 + p_3q_4 + p_5q_5)p_4p_5} \\
&= \frac{q_5}{p_4} \frac{(p_4q_4)p_5 + p_3q_4p_5 + (p_4q_4)p_3 + (p_5q_5)q_4}{(p_5q_5)q_4 + (p_5q_5)p_3 + p_3q_4p_5 + (p_5q_5)p_5} \\
&= \frac{q_5}{p_4} \ \text{(using Condition 3 several times)} \\
\frac{h}{(35) h} &= \frac{p_4p_3}{p_4p_5} = \frac{p_3}{p_5}
\end{align*}

In each case $h/\tau h$ has the divisor we want. $\Box$

\smallskip

\begin{proposition} In the above situation, we can construct functions $p_i$ and $q_i$ satisfying Conditions $1,2,3$. \end{proposition}

{\em Proof:} Define 
\[
p_i(a_0,a_1,a_2,a_3,a_4,a_5) = \sum_{j=0}^5 \left( \frac{r_1^j(r_2-r_i)}{z_2 z_i} + \frac{r_2^j(r_i-r_1)}{z_i z_1} + \frac{r_i^j (r_1-r_2)}{z_1 z_2} \right) a_j.
\]
Define $q_i$ by replacing $z_i$ with $-z_i$ everywhere. We must now show that the intersection divisor of the hyperplane cut out by $p_i$ with $V_{\delta}$ is in fact $E_{\{1,2,i\},\{0,0,0\}}$.

To see this, let $q(x) = \sum_{j=0}^5 a_j x^j$, and suppose that $q(r_m) = \dfrac1{z_m}(tr_m+u)$ for $m = 1,2,i$. Then
\begin{align*}
p_i(a_0,a_1,a_2,a_3,a_4,a_5) &= \frac{q(r_1)(r_2-r_i)}{z_2 z_i} + \frac{q(r_2)(r_i-r_1)}{z_iz_1} + \frac{q(r_i)(r_1-r_2)}{z_1 z_2} \\
&= \frac1{z_1 z_2 z_i} \left( (r_2-r_i)(tr_1 + u) + (r_i-r_1)(tr_2 + u) + (r_1-r_2)(tr_i+u) \right) \\
&= 0.
\end{align*}
Hence the lines whose classes appear in $E_{\{1,2,i\},\{0,0,0\}}$ all lie on the hyperplane cut out by $p_i$. 

Similarly we can show that the intersection divisor of the hyperplane cut out by $q_i$ with $V_{\delta}$ is $E_{\{1,2,i\},\{0,0,1\}}$. Now Conditions 1 and 2 are immediate from the construction of $p_i$ and $q_i$. Note also that $p_iq_i$ actually has coefficients in $k(z_i^2)$, as it is invariant under the transposition of the indices $1,2$ as well as the map $z_i \mapsto -z_i$.

Now we must arrange for Condition 3 to hold by multiplying the $p$'s and $q$'s by suitable constants. For $3 \le i \le 6$, consider the rational function $\dfrac{p_iq_i}{p_6q_6}$ on $\Vbar$. Note that its divisor is $0$, so it is a nonzero constant $u_i$. To evaluate this constant, let $P_v \in V(k_v)$ be a point not in the support of this rational function, and note that 
\[
u_i = \frac{p_i(P_v)q_i(P_v)}{p_6(P_v)q_6(P_v)} = \frac{a_v}{b_v},
\]
where $a_v$ is a norm from $k_v(z_i)$ to $k_v(z_i^2)$ and $b_v$ is a norm from $k_v(z_6)$ to $k_v$. 

\begin{lemma}\label{milne} (\cite{milnecft}, VIII.1.11) Let $k$ be a field of characteristic $\ne 2$. An element $c \in k^*$ is the product of a norm from $k(\sqrt{a})$ and a norm from $k(\sqrt{b})$ if and only if, as an element of $k(\sqrt{ab})$, it is a norm from $k(\sqrt{a},\sqrt{b})$. \end{lemma}

Let $L_i = k(z_i^2)$. For any place $w$ of $L_i$, the expression $u_i = a_v(1/b_v)$ exhibits $u_i \in L_i$ as the product of a norm from $(L_i)_w(z_i)$ and a norm from $(L_i)_w(z_6)$. By Lemma \ref{milne}, we see that $u_i$ is a norm from $(L_i)_w(z_i,z_6)$ to $(L_i)_w(z_iz_6)$. By the Hasse Norm Theorem, it follows that $u_i$ is a norm from $L_i(z_i,z_6)$ to $L_i(z_iz_6)$, say of $d_i \in L_i(z_i,z_6)$. Let $\sigma = (0,0,1,1,1,1) \in G_{\Gamma}$. Then $d_i(\sigma d_i) = u_i$, and if we let $d_j$ be the image of $d_i$ under an element of $H_{96} \cap S_6$ transposing $i$ and $j$, we see that 
\[
\frac{p_iq_i}{p_jq_j} = \frac{d_i \sigma d_i}{d_j \sigma d_j},
\]
and thus if we replace $p_i$ by $p_i/d_i$ and $q_i$ by $q_i/(\sigma d_i)$, we have found functions satisfying Condition (3). $\Box$

\smallskip

Once we have constructed $h$ such that $E-\tau E = \text{div}(h/\tau h)$, the rational function $g$ we want will be a function whose divisor is $(E-(h))+\sigma(E-(h))$, where $\sigma$ is {\em any} element of $H_{96} \setminus H_{48}$. The most convenient $\sigma$ to choose is certainly $(0,0,1,1,1,1) \in G_{\Gamma}$, because $\sigma E = -E$. Thus the rational function $g$ we use is a function whose divisor is $-(h \cdot (0,0,1,1,1,1)h)$. We can drop the negative sign, since we have a quaternion algebra, i.e. $(c,g) = (c,1/g)$. So we can take
\[
g = h \cdot (0,0,1,1,1,1) h \ \text{times a constant}
\]
where we know that we will be able to find a constant so that $g$ is invariant under the action of the Galois group.

\begin{remark} It appears at first glance that $g$ will always be a norm from $k(z_6)$ to $k$, which would make the quaternion algebra we have constructed trivial in $\Br V$. But $h$ is not itself defined over $k(z_6)$, so indeed this algebra is nontrivial in general. \end{remark}

\begin{remark} Here is how we arrange for Condition 3 to hold in practice. We know from the proof above that there are constants $d_i \in k(z_i,z_6)$ (which are conjugate to each other under the action of $S_3$ on the indices $3,4,5$) such that $(p_i/d_i)\sigma(p_i/d_i)$ is independent of $i$, for $3 \le i \le 5$; here $\sigma$ is the element $(0,0,1,1,1,1) \in G_{\Gamma}$ sending $z_i \mapsto -z_i$ for $3 \le i \le 6$. 

We consider the {\em normal form} $N_i$ of $p_iq_i$ with respect to a Gr\"obner basis for the ideal generated by the defining equations of $V_{\delta}$. The normal form is a uniquely determined representative of the equivalence class of $p_iq_i$ modulo this ideal; it is clear from the construction of the normal form that the normal forms of $p_iq_i$ and $p_jq_j$ are permuted just as $p_iq_i$ and $p_jq_j$ are, by the action of $S_3$ on the indices. Fixing a monomial appearing in $N_i$ and calling its coefficient $c_i$, we note that $c_i \in k(z_i)$ should satisfy
\[
\frac{c_i}{c_j} = \frac{p_iq_i}{p_jq_j} = \frac{N_{\sigma}(d_i)}{N_{\sigma}(d_j)}
\]
Note also that 
\[
\frac{N_{\sigma}(d_i)}{c_i}
\]
is independent of $i$, but since this quantity is in $k(z_i^2,z_6)^*$ for each $i$, it follows that it lies in $k(z_6)^*$.

So we search for $n \in k(z_6)^*$ such that $c_i n$ is a norm from $k(z_i,z_6)$ to $k(z_i^2,z_iz_6)$; we are guaranteed that such $n$ exist, and in practice we find them quickly (e.g. using the \texttt{NormEquation} function in MAGMA). 

In fact, in every explicit example the author has written down, he has been able to find $n \in k^*$; but the author does not know if it is always possible to find $n \in k^*$ in general. \end{remark}

\section{Proof of the main theorem}

\begin{theorem} Let $C$ be the hyperelliptic curve $y^2 = (x^2-5x+1)(x^3-7x+10)(x+1)$. Then $\Sha(\Q, \Jac(C))[2] = \Z/2 \times \Z/2$. \end{theorem}

{\em Proof:} Let $J = \Jac(C)$. Let $f(x)$ be the sextic polynomial defining $C$. First note that it is easy to show that $J(\Q)/2J(\Q)$ is a vector space over $\Z/2$ of dimension at least $2$; $J(\Q)[2] = \Z/2$ because there is a quadratic factor of $f(x)$, and there is also the image of a point in $J(\Q)$ coming from the point $(-1,0) \in C(\Q)$. According to Stoll's algorithm from \cite{stoll}, now implemented in MAGMA as \texttt{TwoSelmerGroup}, we can compute $\Sel^{(2)}(\Q,J)$, and find that it is a vector space over $\Z/2$ of dimension equal to $4$. Moreover, since $C$ has a rational point, its Jacobian is even (in the language of \cite{poonenstoll}); that is, the order of $\Sha(\Q,J)[2]$ is a square. Hence it is either $1$ or $4$. We must therefore only exhibit one nontrivial element of $\Sha(\Q,J)[2]$ to complete the proof.

Define $\delta \in A_f^*$ by
\begin{equation}\label{delta}
\delta(x) = -\frac7{2965}(377x^5 - 706x^4 - 5200x^3 + 2061x^2 - 9086x - 12308)
\end{equation}
Then MAGMA shows that $\delta$ gives an element of $\Sel^{(2)}(J)$. The K3 surface $V_{\delta}$ is given by the vanishing of the polynomials
\begin{align*}
& -377a_0^2 - 1604a_0a_1 - 4310a_0a_2 - 9600a_0a_3 - 14130a_0a_4 - 24100a_0a_5 - 2155a_1^2 - 9600a_1a_2 \\
&- 14130a_1a_3 - 24100a_1a_4 + 3002a_1a_5 - 7065a_2^2 - 24100a_2a_3 + 3002a_2a_4 - 3752a_2a_5 \\
&+ 1501a_3^2 - 3752a_3a_4 + 380254a_3a_5 + 190127a_4^2 + 505356a_4a_5 + 2697585a_5^2, \\
& 353a_0^2 + 1053a_0a_1 + 3820a_0a_2 + 12135a_0a_3 + 16210a_0a_4 + 49701a_0a_5 + 1910a_1^2 + 12135a_1a_2 \\
&+ 16210a_1a_3 + 49701a_1a_4 - 7880a_1a_5 + 8105a_2^2 + 49701a_2a_3 - 7880a_2a_4 + 197631a_2a_5 \\
&- 3940a_3^2 + 197631a_3a_4 - 507830a_3a_5 - 253915a_4^2 + 1686873a_4a_5 - 2233480a_5^2, \\
& 1300a_0^2 + 6321a_0a_1 + 17920a_0a_2 + 34505a_0a_3 + 63708a_0a_4 + 62335a_0a_5 + 8960a_1^2 + 34505a_1a_2 \\
&+ 63708a_1a_3 + 62335a_1a_4 + 90560a_1a_5 + 31854a_2^2 + 62335a_2a_3 + 90560a_2a_4 - 243597a_2a_5 \\ 
&+ 45280a_3^2 - 243597a_3a_4 - 202262a_3a_5 - 101131a_4^2 - 3623209a_4a_5 - 1919025a_5^2.
\end{align*}

Since $\delta(r_6) = -7$, we will be using the Azumaya algebra $(-7,g)$, where $g$ is constructed as in the previous section. Note that the denominator of $g$, which is a constant multiple of $h \cdot (0,0,1,1,1,1)h$, is $p_4p_5q_4q_5$. This is a constant multiple of $(p_6q_6)^2$. Plugging in any point in $V(k_v)$ to $(p_6q_6)^2$ yields a norm from $k_v(z_6)$ to $k_v$, so when we compute the invariant of $g$, we can ignore the contribution coming from its denominator. The numerator is a homogeneous polynomial $F$ of degree $4$, so we are reduced now to computing the expressions $(-7,F(P_v))_v$ for all places $v$ of $\Q$ and all points $P_v \in V(k_v)$. (These expressions are invariant under projective scaling of the coordinates of $P_v$, since $F$ has even degree.) 

We wish to compute the primes of bad reduction for $V_{\delta}$ next, for the purposes of invariant computations. This is a standard computation using the minors of the matrix of partial derivatives; we can also note that the only bad primes should be primes dividing the discriminant of the minimal Galois extension over which all the lines are defined. In our case the bad primes are $2,3,7,83,739,\infty$. Note that if $-7$ is a square in $\Q_p$ we know automatically that the quaternion algebra $(-7,F/G)$ will have constant invariant $0$ at $p$. So we do not have to analyze $p=2,739$. 

These will not be the only primes to consider when we compute the possibilities for $\inv_p \calA(P_p)$, $P_p \in V_{\delta}(\Q_p)$. Here we repeat a standard lemma that allows us to restrict our consideration to a finite set of primes.

\begin{lemma} Let $V$ be a smooth projective geometrically integral $k$-variety, $k$ a number field, $g \in k(V)^*$, and let $(L/k,g)$ be a quaternion algebra in $\Br V$. Let $v$ be a nonarchimedian place of $k$ and suppose that $v$ does not ramify in $L$ and that $V$ has smooth reduction at $v$. Then $\inv_v (L/k,g)(P_v)$ is independent of the choice of $P_v \in V(k_v)$. \end{lemma}

{\em Proof:} See Lemma 8.4, \cite{cornlms}. $\Box$

\smallskip

So if $p$ is not a bad prime for $V_{\delta}$ and $p$ does not ramify in $k(z_6)$, then the invariants of the quaternion algebra $(\delta(r_6),g)(P_p)$ will be independent of the choice of $P_p$. When will it be nonzero?

The proof of Lemma 8.4 of \cite{cornlms} implies that, if $\delta(r_6)$ is not a square mod $p$, the invariant of $(\delta(r_6),g)(P_p)$ will equal $m_p/2 \in \frac12\Z/\Z$, where $m_p$ is the integer such that $m_pV_p$ appears in the divisor of $g$ (here $V_p$ denotes the special fiber of the smooth model $\calV$ of $V_{\delta}$ over $\Spec \Z_p$). 

Pick a point $P \in V_{\delta}({\overline \Q})$. Write $g = F/G$, where $F$ and $G$ are as above; they are homogeneous polynomials in six variables with integer coefficients. Since $G$ is a norm from $k(z_6)(V_{\delta})$ to $k(V_{\delta})$, we can ignore it for the purposes of invariant computations. If there is a prime $\mathfrak p$ over $p$ not dividing $F(P)$, then $P_{\mathfrak p} \in V_p({\overline{{\mathbb F}_p}})$ is not in the support of $g$ (or some rational function obtained from $g$ by multiplying the denominator by a norm from $k(z_6)(V_{\delta})$), so the integer $m$ must equal $0$. So the set of primes at which $\inv_p (\delta(r_6),g)(P_p)$ is not always $0$, independent of $P_p \in V_{\delta}(\Q_p)$, is contained in the union of the set of bad primes with the set of primes of $\Z$ dividing $N(F(P))$, where $N$ denotes the absolute norm down to $\Q$.

This is true for all $P$, so after combining results for various $P$ we obtain the set
\[
\{ 5, 61, 347, 739, 3433, 4337, 6833, 663149189, 69804594311 \}
\]
of primes $p$ at which $V_{\delta}$ has smooth reduction but the integer $m_p$ is possibly nonzero. After throwing out the primes for which $-7$ is a square in $\Q_p^*$, we get 
\[
\{ 5, 61, 3433, 663149189 \}.
\]
For each $p$ in this set we find that $m_p$ is odd. This is not hard to check: as the invariant is constant, we need only lift one point not in the support of $F/G$ in $V_{\delta}(\Q_p)$ to high enough precision in order to evaluate $F$ at it. Since the invariant at each of these primes is $1/2$, we obtain 
\begin{align*}
\sum_{v} \inv_v (-7,F/G)(P_v) =& \inv_3 (-7,F/G)(P_3) + \inv_7 (-7,F/G)(P_7) + \inv_{83} (-7,F/G)(P_{83}) \\ 
&+ \inv_{\infty} (-7,F/G)(P_{\infty}).
\end{align*}

To carry out the invariant computation at $\infty$, we first show that the value of the function $F$ on the points in $V_{\delta}({\mathbb R})$ is either always positive or always negative. To see this, consider the polynomial whose roots are $\pm z_3, \pm z_4, \pm z_5$; this polynomial can be written as $h(x^2)$, where $h(x) = x^3-126x^2+7938x+250047$. Notice that $h(x)$ has one real root, which is negative. Under any embedding of the splitting field of $h(x^2)$ into $\C$ such that the image of $z_3^2$ is the negative real root, the images of $z_4^2$ and $z_5^2$ are distinct and conjugate. Choose such an embedding so that $z_4$ and $-z_5$ are complex conjugates. Then consider the action of complex conjugation on the $p_i$ and $q_i$; we see that it sends $p_3$ to $q_3$ and vice versa, while it sends $p_4$ to $q_5$ and $p_5$ to $q_4$.

Now $F$ is a constant multiple of 
\[
(p_3p_5 + p_3p_4 + p_4p_5 + p_3q_3)(q_3q_5 + q_3q_4 + q_4q_5 + p_3q_3),
\]
and if we plug in a point in $V_{\delta}({\mathbb R})$ to both these factors, we see immediately that we get the product of a complex number and its conjugate. So if we can show that $F$ is never zero on $V_{\delta}({\mathbb R})$, we can deduce that it is either always positive or always negative. It can be easily verified--either by plugging in a specific real-valued point or noticing that $F$ is actually a {\em positive} constant multiple of the above function--that $F$ is always positive. This shows that the invariant is actually zero.

To see that $F$ never vanishes on $V_{\delta}({\mathbb R})$, note that if $F(P) = 0$, then both factors of $F$ vanish at $P$. So $P$ lies on the intersection of the K3 surface $V_{\delta}$ with the two hyperplanes cut out by the factors of $F$; this intersection is zero-dimensional, and MAGMA computes that there are no real points on it.

At $83$ we note that there is one singular ${\mathbb F}_{83}$-valued point in the special fiber $V_{83}$, but it does not lift to any points mod $83^2$. By a theorem of Bright (\cite{brighteff}, Theorem 1), the invariant is constant above any smooth ${\mathbb F}_{83}$-point, so we need merely write down all points in $V_{\delta}({\mathbb F}_{83})$, and lift each one to high enough precision in order to evaluate $F/G$ at it. (Hensel's lemma guarantees that any lift of a smooth ${\mathbb F}_{83}$-point will lift to a $\Z_{83}$-point.) Doing this computation for each of the $6960$ smooth points in $V_{\delta}({\mathbb F}_{83})$, we find that each point has a lift $P$ mod $83^5$ such that the $83$-adic valuation of $F(P)$ is either $2$ or $4$. The invariant at $P$ is $0$ or $1/2$ depending on whether or not $F(P)$ has even or odd valuation, respectively; so in all cases, we find that the invariant at $83$ is zero.

At $7$, we compute that $V({\mathbb F}_7)$ has $71$ points. One of these points, $P_s = (3 \colon 6 \colon 3 \colon 1 \colon 2 \colon 1)$, is singular, and the set 
\[ 
\{ P \in V({\mathbb F}_7) \colon F(P) = 0 \}
\]
has $15$ elements (including the singular point). For the other $56$ points, we find that $F(P) = 1$, $2$, or $4$, which are the squares in ${\mathbb F}_7^*$, so the invariant is zero at all these points. 

Again using Theorem 1 of \cite{brighteff}, we find one lift $P$ of each of the $14$ remaining nonsingular points in $V({\mathbb F}_7)$ to a point in $V(\Z_7)$ (at least one of whose coordinates is a unit--we will call such $P$ {\em normalized}) and compute $F(P)$ for each such lift. In particular, we see that $F(P) \equiv 7^2$, $2 \cdot 7^2$, or $4 \cdot 7^2$ mod $7^3$ for each of these points, so $F(P)$ is a square in $\Z_7$, so the invariant is zero at these points as well.

Above the point $P_s$, Bright's theorem does not apply, and we must actually consider all possible lifts of $P_s$ to $P \in V(\Z_7)$. A lengthy MAGMA computation shows that for any such normalized $P$, $F(P) \equiv 2 \cdot 7^4$ mod $7^5$, so once again the invariant is zero at all these points, and hence the invariant is zero at all points in $V(\Q_7)$.

At $3$, we compute that $V({\mathbb F}_7)$ has $40$ points. All of them are singular, so Bright's theorem does not apply. A MAGMA computation similar to the one at $7$ shows that if $P$ is any normalized point in $V(\Z_3)$, $v_3(F(P)) = 5$. So the invariant is $1/2$ at all points in $V(\Q_3)$. 

Hence the sum $\sum \inv_v \calA(P_v)$ is constant and equal to $1/2$ for all $(P_v) \in V_{\delta}({\mathbb A}_\Q)$, so $V_{\delta}({\mathbb A}_\Q)^{\rm Br} = \emptyset$, so $V_{\delta}(\Q) = \emptyset$. So $\delta \in \Sel^{(2)}(\Q,J)$ maps to a nonzero element of $\Sha(\Q,J)[2]$. $\Box$

\begin{proposition}\label{mainthm} Let $n$ be a (positive) product of primes splitting completely in the field 
\[
L(\sqrt{2},\sqrt{-739}),
\]
where $L$ is the field of definition of the $32$ lines, the degree-$96$ splitting field of $(x^2-35x+49)(x^6-126x^4+7938x^2+250047)$. Then the Jacobian $J_n$ of the curve given by
\[
y^2 = n(x^2-5x+1)(x^3-7x+10)(x+1)
\]
satisfies $\Sha(\Q,J_n)[2] \ne 0$. \end{proposition}

{\em Proof:} Let $f(x) = (x^2-5x+1)(x^3-7x+10)(x+1)$. Note that $A_f = A_{nf}$. We show that the element 
\[
\delta = -\frac{7}{2965}(377x^5 - 706x^4 - 5200x^3 + 2061x^2 - 9086x - 12308)
\]
lies in $\Sel^{(2)}(\Q,J_n)$. By Corollary 5.11 of \cite{stoll}, we must check this only at $\infty$, $2$, and primes $p$ such that $p^2$ divides the discriminant of $nf(x)$. So we must check $\infty$, $2$, $7$, $739$, and primes dividing $n$. At primes $p$ dividing $n$, we actually show that the image of $\delta$ lies in $(A_{nf}^*)^2 \Q_p^*$. Let $r_1$ and $r_2$ be the roots of $x^2-5x+1$, let $r_3,r_4,r_5$ be roots of $x^3-7x+10$, and let $r_6 = -1$. Consider the isomorphism
\[
A_{nf} \to \Q(r_1) \oplus \Q(r_3) \oplus \Q 
\]
given by $g(x) \mapsto (g(r_1), g(r_3), g(r_6))$. The image of $\delta$ under this isomorphism is 
\[
(\delta(r_1),\delta(r_3),\delta(-1)) = (z_1^2,z_3^2,-7).
\]
For primes $p$ splitting completely in $L$, these three elements $z_1^2$, $z_3^2$, and $-7$ are actually in $(\Q_p^*)^2$. So in fact the image of $g(x)$ lies in $((A_{nf} \otimes \Q_p)^*)^2$, so it is zero in $H_f \otimes \Q_p$; hence it is trivially in the image of $\Delta_{\Q_p}$. 

At the primes $p = \infty, 2, 7, 739$, we notice that $n$ is a square in $\Q_p^*$; so the fact that $\delta$ lies in $\Delta_{\Q_p}(J_n(\Q_p)/2J_n(\Q_p))$ follows from the fact that it lies in $\Delta_{\Q_p}(J_1(\Q_p)/2J_1(\Q_p))$. This completes the proof. $\Box$

\section{Connections with Del Pezzo surfaces of degree $4$}

In this section, we consider the odd-degree model $C'$ of $C$ obtained by sending the rational point $(-1,0)$ to infinity. Its equation is
\[
w^2 = (7u^2-7u+1)(16u^3-4u^2-3u+1),
\]
via the change of variables $w=y/(x+1)^3$, $u=1/(x+1)$. Here we prove a general lemma about the relationship between these two curves.

\begin{lemma} Consider a hyperelliptic curve $C_1$ over a number field $k$ with equation
\[
y^2 = f_5(x)(x-a).
\]
Let $C_2$ be the odd-degree model of this curve, given by the equation
\[
w^2 = u^5 f_5\left( \frac1{u}+a \right)
\]
where $w = y/(x-a)^3$, $u=1/(x-a)$. Let $f(x)$ be the sextic polynomial defining $C_1$ and $g(u)$ the quintic polynomial defining $C_2$. Let $\delta \in \Sel^{(2)}(k,\Jac C_1)$, considered as in \cite{stoll} as a subgroup of $A_f^*/(A_f^*)^2 k^*$, and let $\beta$ be an element in $\Sel^{(2)}(k,\Jac C_2)$, considered as in \cite{stoll} as a subgroup of $A_g^*/(A_g^*)^2$. There is a natural isomorphism $\Sel^{(2)}(k,\Jac C_1) \to \Sel^{(2)}(k,\Jac C_2)$; suppose it sends $\delta$ to $\beta$. Let $W_{\beta}$ be the intersection of two quadrics obtained as the projectivization of the set
\[
\calW_{g,\beta} = \{ s(u) \in A_g \colon \beta(u)s(u)^2 \equiv \text{quadratic (mod $g$)} \}
\]
(cf. \cite{bruinflynn} and \cite{logandp4}). Then the K3 surface $V_{\delta}$ obtained from $C_1$ is a double cover of $W_{\beta}$, which is a smooth Del Pezzo surface of degree $4$. It ramifies over the projectivization of
\[
\{ s(u) \in A_g \colon \beta(u)s(u)^2 \equiv \text{linear (mod $g$)} \}.
\]
\end{lemma}

{\em Proof:} Note that $k[x]/(f_5(x)) \cong k[u]/(g(u))$. Now there is a natural isomorphism between $A_f^*/(A_f^*)^2 \Q^*$ and $A_g/(A_g^*)^2$, as the first group is isomorphic to 
\[
\frac{(k[x]/f_5(x))^* \oplus k^*}{\text{squares} \cdot k^*}
\]
and the second group is isomorphic to 
\[
\frac{(k[u]/(g(u)))^*}{\text{squares}}
\]

It is not hard to see that this isomorphism restricts to an isomorphism $\phi \colon \Sel^{(2)}(k, \Jac C_1) \to \Sel^{(2)}(k,\Jac C_2)$ (where these Selmer groups are realized, as in \cite{stoll}, as subgroups of $A_f^*/(A_f^*)^2 k^*$ and $A_g/(A_g^*)^2$, respectively). 

The fact that $W_{\beta}$ is a smooth Del Pezzo surface of degree $4$ is Lemma 17 of \cite{logandp4}. Now suppose $q \in \calV_{f,\delta}$, so that $\delta q^2$ is congruent to a quadratic polynomial $c$ mod $f$. Let $c_1(u) = u^2 c(1/u+a)/\delta(a)$; then $c_1$ is also quadratic. If $\beta = \phi(\delta)$, then
\[
\beta(u) \equiv u^6 \delta(1/u+a)/\delta(a) \ \text{mod $(A_g^*)^2$},
\]
and if we choose this representative polynomial for the class of $\beta$, we see that the polynomial $s(u)$ defined by the formula
\begin{equation}\label{doublecover}
s(u) \equiv u^{-2} q(1/u+a) \ \text{mod $g(u)$}
\end{equation}
is in $\calW_{g,\beta}$, because $\beta(u) s(u)^2 \equiv c_1(u)$ (mod $g(u)$). Note that $s(u)$ is well-defined because $u$ is invertible mod $g(u)$. 

So (\ref{doublecover}) gives a map $\calV_{f,\delta} \to \calW_{g,\beta}$, which induces a map $V_{\delta} \to W_{\beta}$. Now consider $s(u) \in \calW_{g,\beta}$, and suppose that $\beta s^2 \equiv c_1$ mod $g$. Let $r_1, \ldots, r_6$ be the roots of $f$ in $\kbar$, where $r_6=a$. Now $q \in \calV_{f,\delta}$ maps to $s$ if and only if
\[
q(r_i) = \frac1{(r_i-a)^2} s\left( \frac1{r_i-a} \right), \quad 1 \le i \le 5.
\]
Note that the quadratic polynomial associated to $q$ is $c(x) = \delta(a) (x-a)^2 c_1(1/(x-a))$. So if $q$ is to lie in $\calV_{f,\delta}$, we must have that $q(a)^2$ is the leading term of $c_1$. Generically there are two choices for the square root of this leading term (these choices coincide when the leading term is $0$). The result follows. $\Box$

\smallskip

As mentioned above, the lemma implies that if $W_{\beta}$ has a Brauer-Manin obstruction to rational points, then so does $V_{\delta}$, and if $W_{\beta}$ fails the Hasse principle, then so does $V_{\delta}$. So the ``K3 method" of exhibiting explicit elements of $\Sha[2]$ given above can be viewed as a generalization of the ``Del Pezzo method" of \cite{bruinflynn} and \cite{logandp4}. (Of course, the K3 method also works more generally; it applies to any genus-$2$ curve $y^2=f(x)$, while the Del Pezzo method works only for curves with odd-degree models. The example given in \cite{ronaldadam} is an example of an application of the K3 method to a curve without an odd-degree model.)

For the curve $y^2 = (x^2-5x+1)(x^3-7x+10)(x+1)$ and specific choice of $\delta$ given in (\ref{delta}), this construction produces a Del Pezzo surface $W_{\beta}$ with points everywhere locally. As it happens, however, Logan's program for Del Pezzo surfaces (\cite{logan}) shows that $W_{\beta}$ has no Brauer-Manin obstruction to rational points. We expect that it satisfies the Hasse principle as well (although the size of the coefficients in the defining equations for $W_{\beta}$ precludes a successful search for rational points of small naive height). In any case, this shows that the K3 method does in fact give a generalization of the Del Pezzo method even in the case when both methods apply. 

\bibliographystyle{plain}
\bibliography{pat}

\end{document}